\documentclass[oneside,english]{amsart}
\usepackage[T1]{fontenc}
\usepackage[latin9]{inputenc}
\usepackage{textcomp}
\usepackage{amsthm}
\usepackage{amstext}
\usepackage{amssymb}

\makeatletter
\numberwithin{equation}{section}
\numberwithin{figure}{section}
\theoremstyle{plain}
\newtheorem{thm}{\protect\theoremname}
  \theoremstyle{definition}
  \newtheorem{problem}[thm]{\protect\problemname}
  \theoremstyle{definition}
  \newtheorem{defn}[thm]{\protect\definitionname}

\makeatother

\usepackage{babel}
  \providecommand{\definitionname}{Definition}
  \providecommand{\problemname}{Problem}
\providecommand{\theoremname}{Theorem}

\begin{document}

\title{Calabi-Yau structures on cotangent bundles}

\author{Alexandru Doicu}
\begin{abstract}
Starting with a orientable compact real-analytic Riemannian manifold
$(L,g)$ with $\chi(L)=0$, we show that a small neighbourhood $\textnormal{Op}(L)$
of the zero section in the cotangent bundle $T^{*}L$ carries a Calabi-Yau
structure such that the zero section is an isometrically embedded
special Lagrangian submanifold.
\end{abstract}
\maketitle
\tableofcontents{}

\section{\label{sec:Introduction}Introduction}

Let $(L,g)$ be an orientable compact real-analytic Riemannian manifold
with real-analytic Riemannian metric $g$. According to \cite{key-1}
there exists a sufficiently small neighbourhood $\textnormal{Op}(L)$
of the zero section (which we identify with $L$) in the cotangent
bundle $T^{*}L$ carrying a complex structure $J$ which we view as
a integrable almost complex structure. With respect to this complex
structure, $L\hookrightarrow\textnormal{Op}(L)$ is a totally real
submanifold. On $\textnormal{Op}(L)$ there exists a strictly plurisubharmonic
exhaustion function $\rho:\textnormal{Op}(L)\rightarrow\mathbb{R}$
such that the Kähler metric $\tilde{g}$ obtained from the Kähler
form $\omega=(i/2)\partial\overline{\partial}\rho$ restricts to the
Riemannian metric $g$ when restricted to the zero section $L$ \cite{key-2}.
Summarizing, we are dealing with a Kähler manifold $(\textnormal{Op}(L),J,\omega)$
with Kähler metric $\tilde{g}$, such that $\omega|_{L}\equiv0$ and
$\tilde{g}|_{L}\equiv g$.

We ask if it is possible to solve the following problem:
\begin{problem}
Find a pair $(J,\omega)$ and a non-vanishing holomorphic $(n,0)$-form
$\Omega$ on $\textnormal{Op}(L)$ with the properties:
\begin{enumerate}
\item $(\textnormal{Op}(L),J,\omega)$ is a Kähler manifold.
\item $\omega|_{L}\equiv0,\tilde{g}|_{L}\equiv g$ and $\Omega|_{L}\equiv\textnormal{Vol}_{g}$,
where $\textnormal{Vol}_{g}$ is the volume form on $L$ induced by
the Riemannian metric $g$, and $\tilde{g}$ is the Kähler metric
$\tilde{g}(\cdot,\cdot)=\omega(\cdot,J\cdot)$.
\item $\frac{\omega^{n}}{n!}=(-1)^{\frac{n(n-1)}{2}}\left(\frac{i}{2}\right)^{n}\Omega\wedge\overline{\Omega}$.
\end{enumerate}
\end{problem}
Some of these structures are already uniquely determined. For example,
the complex structure $J$ is unique on a suffciently small neigbourhood
of the $0$-section in the cotangent bundle of $L$ up to biholomorphism
\cite{key-1}. Moreover, on a suffciently small neigbourhood $\textnormal{Op}(L)$,
there exists a unique holomorphic $(n,0)$-form with the property
$\Omega|_{L}=\textnormal{Vol}_{g}$. Indeed, if the Riemannian metric
$g$ on $L$ is real-analytic, then the volume form $\textnormal{Vol}_{g}$
is given locally by the real-analytic function $\sqrt{\det(g_{ij})}$.
This function can be holomorphically extended, and by the compactness
of $L$, this extension process gives rise to a holomorphic $(n,0)$-form
$\Omega$ such that $\Omega|_{L}=\textnormal{Vol}_{g}$. Hence, the
complex structure $J$ together with the holomorphic $(n,0)$-form
$\Omega$ are unique. However, the Kähler form $\omega$ is not necessarily
unique, and Problem 1 reduces to finding a specific Kähler form $\omega$
satisfying the requirements (1), (2) and (3). A Kähler manifold which
admits a holomorphic $(n,0)$-form satisfying Equation (3) is called
a \emph{Calabi-Yau }manifold. If, in addition, the conditions $\omega|_{L}\equiv0$
and $\Omega|_{L}\equiv\textnormal{Vol}_{g}$ are fulfilled, the submanifold
$L\hookrightarrow\textnormal{Op}(L)$ will be called a \emph{special
Lagrangian }submanifold. Problem 1 has up to now resisted a complete
solution. However, it has been solved for some special classes of
manifolds:
\begin{enumerate}
\item For a three-dimensional, compact, real-analytic Riemannian manifold,
Bryant \cite{key-3} solved Problem 1 by using techniques from Cartan-Kähler
theory. The main idea in his proof is that every compact, real-analytic
Riemannian manifold with real-analytic metric is real-analytically
parallelizable, and that Problem 1 can be reduced to a problem of
finding particular integral submanifolds of an exterior differential
ideal.
\item For a compact, real-analytic Kähler manifold with real-analytic Kähler
form and complex structure, Feix \cite{key-4} proved that a neighbourhood
of the $0$-section of its cotangent bundle carries a Hyperkähler
structure, and demonstrated that by rotating the complex structures
together with the Kähler forms we are led to a Calabi-Yau structure
on a neigbourhood of the cotangent bundle such that the $0$-section
is a special Lagrangian submanifold.
\item For a compact, rank one, globally symmetric space $L$, Stenzel \cite{key-5}
showed that Equation (3) can be reduced, by using symmetries, to a
solvable ordinary differential equation, and that $L$ is an isometrically
embedded special Lagrangian submanifold.
\end{enumerate}
The main result of this paper is the solution of Problem 1 in the
sace $\chi(L)=0$ (generalising Bryant's result):
\begin{thm}
\label{thm:Let--be--MAIN-ONCE-AGAIN}Let $(L,g)$ be an orientable
compact real-analytic Riemannian manifold with $\chi(L)=0$. Then,
there exists a Calabi-Yau structure $(\textnormal{Op}(L),\Omega,\omega)$,
where $\textnormal{Op}(L)\subset T^{*}L$, such that $L$ is a isometrically
embedded special Lagrangian submanifold.
\end{thm}
The paper is organized as follows: In Section \ref{sec:Calabi-yau-manifolds-and},
we recall some geometric properties of Calabi-Yau manifolds and special
Lagrangian submanifolds. In Section \ref{sec:Proof-of-Theorem}, we
prove Theorem \ref{thm:Let--be--MAIN-ONCE-AGAIN}.

\section{\label{sec:Calabi-yau-manifolds-and}Calabi-yau manifolds and special
Lagrangian submanifolds}

Calabi-Yau manifolds can be defined in several ways. First, they can
be regarded as Kähler manifolds equipped with a Ricci-flat Kähler
metric. In the compact case, we can drop the assumption on the Ricci-flatness
and define them to have vanishing first Chern class \cite{key-8}.
Second, they can be defined as Kähler manifolds with holonomy group
contained in $\textnormal{SU}(n)$. Third, they can be introduced
via the following: 
\begin{defn}
\label{def:A-K=0000E4hler-manifold}A Kähler manifold $(X,J,\omega)$
is called \emph{Calabi-Yau} if there exists a nonvanishing holomorphic
$(n,0)$-form $\Omega$ such that the following equation holds
\begin{equation}
\frac{\omega^{n}}{n!}=(-1)^{\frac{n(n-1)}{2}}\left(\frac{i}{2}\right)^{n}\Omega\wedge\overline{\Omega}.\label{eq:1}
\end{equation}

\end{defn}
This definition will be adopted in our analysis. Some properties of
Calabi-Yau manifolds are summarized below.
\begin{enumerate}
\item \label{enu:The-constant-}The constant $(-1)^{n(n-1)/2}(i/2)^{n}$
is a normalisation constant, so that in local holomorphic coordinates
$(z_{1},...,z_{n})$, (\ref{eq:1}) is of the form 
\begin{equation}
\det\left(g_{i\overline{j}}\right)=|h|^{2},\label{eq:2}
\end{equation}
where $g_{i\overline{j}}$ are the coefficients of the Kähler form
$\omega$ in the basis $\partial/\partial z^{i},\partial/\partial\overline{z}^{j}$
and $h$ is a holomorphic function such that $\Omega=hdz^{1}\wedge...\wedge dz^{n}$.
\item \label{enu:A-Calabi-Yau-manifold}A Calabi-Yau manifold is Ricci-flat.
This can be seen as follows. The Ricci-form in local holomorphic coordinates
is given by $\textnormal{Ric}(\omega)=-i\partial\overline{\partial}\textnormal{log}(\det(g_{i\overline{j}}))$.
Thus, from (\ref{eq:2}) we obtain 
\begin{align*}
\textnormal{Ric}(\omega) & =-i\partial\overline{\partial}\log(\det(g_{i\overline{j}}))=-i\partial\overline{\partial}\log(|h|^{2})\\
 & =-i\partial\overline{\partial}\log(h\overline{h})=0.
\end{align*}

\item \label{enu:The-holonomy-group}The holonomy group of a Calabi-Yau
manifold is contained in $\textnormal{SU}(n)$. To see this, note
that a Calabi-Yau manifold has trivial canonical bundle $K_{X}=\Lambda^{(n,0)}T^{*}X$.
Denote by $g$ the Kähler metric induced by the Kähler form $\omega$.
Furthermore, let $\nabla$ be the Levi-Civita connection of the metric
$g$. Thus, we can induce a metric on the cotangent bundle $T^{*}X$
and in particular, a Riemannian metric on the canonical bundle $K_{X}$.
Denote by $\nabla^{K}$ the Levi-Civita connection on $K_{X}$ induced
by the Riemannian metric. By straightforward calculation we can show
that the Ricci-form $\textnormal{Ric}(\omega)$ is equal to $-i$
times the curvature tensor of the canonical line bundle $K_{X}$.
Thus, the Ricci-form is zero if and only if there exists a parallel,
hence holomorphic, form of type $(n,0)$ in a neighbourhood of any
point of $X$. For Calabi-Yau manifolds this is obviously true. Consequently,
from $\nabla g=0,\nabla\omega=0$ and $\nabla\Omega=0$ we infer that
the holonomy on the Levi-Civita connection induced by the Riemannian
metric $g$ is contained in $\textnormal{SU}(n)$. For detailed calculations
the reader might consult \cite{key-7}. 
\end{enumerate}
Special Lagrangian submanifolds of Calabi-Yau manifolds were first
introduced by Harvey and Lawson in \cite{key-9} as a particular case
of calibrated submanifolds. There are some equivalent definitions
of special Lagrangian submanifolds. In our analysis we will use the
following:
\begin{defn}
\label{def:A-submanifold-}A submanifold $Y$ of a Calabi-Yau manifold
$\left(X,J,\omega,\Omega\right)$ is \emph{special Lagrangian} if\emph{
\[
\omega|_{L}=0\text{ and }\textnormal{Im}(\Omega)|_{L}=0.
\]
}
\end{defn}

\section{\label{sec:Proof-of-Theorem}Proof of Theorem \ref{thm:Let--be--MAIN-ONCE-AGAIN}}

To prove Theorem \ref{thm:Let--be--MAIN-ONCE-AGAIN}, we use the fact
that every compact manifold $L$ with $\chi(L)=0$ admits a nonvanishing
vector field and conversely \cite{key-10}.

We come now to the proof. Let $(L,g)$ be a real-analytic Riemannian
manifold with real-analytic metric $g$. Since $\chi(L)=0$, there
exists a globally defined non-vanishing vector field $X$ on $L$.
According to \cite{key-11}, we can choose $X$ to be real analytic.
Let $V_{x}=\textnormal{span}\{X(x)\}$ and $W_{x}=\{Y\in T_{x}L|g_{x}(Y,X(x))=0\}$,
and let $V=\coprod_{x\in L}V_{x}$ and $W=\coprod_{x\in L}W_{x}$.
$V$ and $W$ are real-analytic vector bundles over $L$ of rank $1$
and $n-1$, respectively, and it is apparent that $V\oplus W=TL$.
If we regard $V$, $W$ and $V\oplus W$ as manifolds, we have that
$\dim(V)=n+1$, $\dim(W)=2n-1$ and $\dim(V\oplus W)=2n$. Let us
consider the map 
\[
TL=V\oplus W\overset{F}{\longrightarrow}W\times\mathbb{R},
\]
given by $F(p,tX(p),Y)=(p,Y,t)$. This map is a real-analytic diffeomorphism.
$L$ is identified with the zero section $L\times\{0\}$ in $TL$,
and via $F$, $L$ is identified with $L\times\{0\}\times\{0\}\subset W\times\mathbb{R}$.
Consider the Kähler structure introduced by Stenzel \cite{key-2}
on $TL$, i.e., $J:T\left(TL\right)\rightarrow T\left(TL\right)$
with the Kähler form $\left(i/2\right)\partial\bar{\partial}\rho$,
where $\rho$ is a strictly plurisubharmonic exhaustion function defined
on $\textnormal{Op}(L)$. By $F$, we transport these structures to
$W\times\mathbb{R}$, so that $F$ becomes an isometric biholomorphism.
To avoid an abundance of notations, denote again by $J$ the complex
structure on $W\times\mathbb{R}$, and by $\rho$ the strictly plurisubharmonic
function which defines the Kähler form. Consider now the Kähler manifold
$(W\times\mathbb{R},J,\left(i/2\right)\partial\overline{\partial}\rho)$,
where the metric induced by $\left(i/2\right)\partial\overline{\partial}\rho$
restricts to $g$ on $L\times\{0\}\times\{0\}$. We intend to find
a plurisubharmonic function $\phi$ in a neighbourhood of $L\times\{0\}\times\{0\}$
in $W\times\mathbb{R}$, such that the Kähler metric induced by $\left(i/2\right)\partial\overline{\partial}\phi$
restricts on $L$ to $g$, and such that 
\begin{equation}
\left(\frac{i}{2}\partial\overline{\partial}\phi\right)^{n}=(-1)^{\frac{n(n-1)}{2}}\left(\frac{i}{2}\right)^{n}\Omega\wedge\overline{\Omega}.\label{eq:AP1}
\end{equation}

In order to apply the Cauchy-Kovalewsky theorem, we equip (\ref{eq:AP1})
with the initial conditions $\phi(p,0)=\rho(p,0)$ and $\frac{\partial\phi}{\partial t}(p,0)=\frac{\partial\rho}{\partial t}(p,0)$.
Hence we are looking for solutions to the following initial value
problem
\begin{equation}
\begin{cases}
\left(\frac{i}{2}\partial\overline{\partial}\phi\right)^{n} & =(-1)^{\frac{n(n-1)}{2}}\left(\frac{i}{2}\right)^{n}\Omega\wedge\overline{\Omega}\\
\phi(p,0) & =\rho(p,0)\\
\frac{\partial\phi}{\partial t}(p,0) & =\frac{\partial\rho}{\partial t}(p,0).
\end{cases}\label{eq:first-pde}
\end{equation}

For $p_{0}\in L$, let $(x^{1},...,x^{n})$ be the normal coordinates
in a neighbourhood $U$ of $p_{0}$ in $L$. Moreover, let $(Y_{1},...,Y_{n-1},X)$
be an orthonormal frame in $U$, with respect to $g$, such that
\[
(Y_{1}(p_{0}),...,Y_{n-1}(p_{0}),X(p_{0}))=(\frac{\partial}{\partial x^{1}}|_{p_{0}},...,\frac{\partial}{\partial x^{n-1}}|_{p_{0}},\frac{\partial}{\partial x^{n}}|_{p_{0}}).
\]
Hence,
\begin{equation}
(x^{1},...,x^{n},y^{1},...,y^{n-1},t)\mapsto(x^{1},...,x^{n},\sum_{i=1}^{n-1}y^{i}Y_{i}+tX)\label{eq:coordinates-choice}
\end{equation}
is a coordinate chart around $p_{0}$ in $TL$ with the property that
\[
(x^{1},...,x^{n},y^{1},...,y^{n-1})\mapsto(x^{1},...,x^{n},\sum_{i=1}^{n-1}y^{i}Y_{i})
\]
is a trivialisation of the bundle $W\rightarrow L$. The coordinates
(\ref{eq:coordinates-choice}) are real-analytic on $W\times\mathbb{R}$,
where $t$ is the global coordinate on $\mathbb{R}$. In these coordinates,
we have 
\begin{align}
J_{p_{0}}\frac{\partial}{\partial x^{i}}|_{p_{0}} & =\frac{\partial}{\partial y^{i}}|_{p_{0}}\text{ for }i=1,...,n-1\label{eq:complex-structure-in-p}\\
J_{p_{0}}\frac{\partial}{\partial x^{n}}|_{p_{0}} & =\frac{\partial}{\partial t}|_{p_{0}}\label{eq:complex-structure-in-p-1}
\end{align}
and the complex structure takes the form 
\[
J=\left(\begin{array}{cccccccc}
J_{x^{1}}^{x^{1}} & J_{x^{2}}^{x^{1}} & ... & J_{x^{n}}^{x^{1}} & J_{y^{1}}^{x^{1}} & ... & J_{y^{n-1}}^{x^{1}} & J_{t}^{x^{1}}\\
J_{x^{1}}^{x^{2}} & J_{x^{2}}^{x^{2}} & ... & J_{x^{n}}^{x^{2}} & J_{y^{1}}^{x^{2}} & ... & J_{y^{n-1}}^{x^{2}} & J_{t}^{x^{2}}\\
\vdots & \vdots & \ddots & \vdots & \vdots & \ddots & \vdots & \vdots\\
J_{x^{1}}^{x^{n}} & J_{x^{2}}^{x^{n}} & ... & J_{x^{n}}^{x^{n}} & J_{y^{1}}^{x^{n}} & ... & J_{y^{n-1}}^{x^{n}} & J_{t}^{x^{n}}\\
J_{x^{1}}^{y^{1}} & J_{x^{2}}^{y^{1}} & ... & J_{x^{n}}^{y^{1}} & J_{y^{1}}^{y^{1}} & ... & J_{y^{n-1}}^{y^{1}} & J_{t}^{y^{1}}\\
\vdots & \vdots & \ddots & \vdots & \vdots & \ddots & \vdots & \vdots\\
J_{x^{1}}^{y^{n-1}} & J_{x^{2}}^{y^{n-1}} & ... & J_{x^{n}}^{y^{n-1}} & J_{y^{1}}^{y^{n-1}} & ... & J_{y^{n-1}}^{y^{n-1}} & J_{t}^{y^{n-1}}\\
J_{x^{1}}^{t} & J_{x^{2}}^{t} & ... & J_{x^{n}}^{t} & J_{y^{1}}^{t} & ... & J_{y^{n-1}}^{t} & J_{t}^{t}
\end{array}\right).
\]
To compute $\left(i/2\right)\partial\overline{\partial}\phi$, we
note that 
\begin{align*}
\frac{\partial}{\partial z^{i}} & =\frac{1}{2}\left(\frac{\partial}{\partial x^{i}}-iJ\frac{\partial}{\partial x^{i}}\right)=\frac{1}{2}\left(\frac{\partial}{\partial x^{i}}-i\left(J_{x^{i}}^{x^{k}}\frac{\partial}{\partial x^{k}}+J_{x^{i}}^{y^{k}}\frac{\partial}{\partial y^{k}}+J_{x^{i}}^{t}\frac{\partial}{\partial t}\right)\right)
\end{align*}
and that 
\[
\frac{\partial}{\partial\overline{z}^{i}}=\frac{1}{2}\left(\frac{\partial}{\partial x^{i}}+iJ\frac{\partial}{\partial x^{i}}\right)=\frac{1}{2}\left(\frac{\partial}{\partial x^{i}}+i\left(J_{x^{i}}^{x^{k}}\frac{\partial}{\partial x^{k}}+J_{x^{i}}^{y^{k}}\frac{\partial}{\partial y^{k}}+J_{x^{i}}^{t}\frac{\partial}{\partial t}\right)\right),
\]
for $i=1,...,n$. Then, we obtain
\begin{align*}
\frac{\partial^{2}\phi}{\partial\overline{z}^{j}\partial z^{i}} & =\frac{1}{2}\frac{\partial}{\partial\overline{z}^{j}}\left(\frac{\partial\phi}{\partial x^{i}}-iJ_{x^{i}}^{x^{k}}\frac{\partial\phi}{\partial x^{k}}-iJ_{x^{i}}^{y^{k}}\frac{\partial\phi}{\partial y^{k}}-iJ_{x^{i}}^{t}\frac{\partial\phi}{\partial t}\right)\\
 & =\frac{1}{4}\left(\frac{\partial}{\partial x^{j}}+iJ_{x^{j}}^{x^{l}}\frac{\partial}{\partial x^{l}}+iJ_{x^{j}}^{y^{l}}\frac{\partial}{\partial y^{l}}+iJ_{x^{j}}^{t}\frac{\partial}{\partial t}\right)\\
 & \times\left(\frac{\partial\phi}{\partial x^{i}}-iJ_{x^{i}}^{x^{k}}\frac{\partial\phi}{\partial x^{k}}-iJ_{x^{i}}^{y^{k}}\frac{\partial\phi}{\partial y^{k}}-iJ_{x^{i}}^{t}\frac{\partial\phi}{\partial t}\right)\\
 & =\frac{1}{4}\left(F_{i,j}+J_{x^{i}}^{t}J_{x^{j}}^{t}\frac{\partial^{2}\phi}{\partial t^{2}}\right),
\end{align*}
where 
\begin{align*}
F_{i,j} & =\frac{\partial^{2}\phi}{\partial x^{j}\partial x^{i}}-i\partial_{x^{j}}J_{x^{i}}^{x^{k}}\frac{\partial\phi}{\partial x^{k}}-iJ_{x^{i}}^{x^{k}}\frac{\partial^{2}\phi}{\partial x^{j}\partial x^{k}}-i\partial_{x^{j}}J_{x^{i}}^{y^{k}}\frac{\partial\phi}{\partial y^{k}}-iJ_{x^{i}}^{y^{k}}\frac{\partial^{2}\phi}{\partial x^{j}\partial y^{k}}\\
 & -i\partial_{x^{j}}J_{x^{i}}^{t}\frac{\partial\phi}{\partial t}-iJ_{x^{i}}^{t}\frac{\partial^{2}\phi}{\partial x^{j}\partial t}+iJ_{x^{j}}^{x^{l}}\frac{\partial^{2}\phi}{\partial x^{l}\partial x^{i}}+J_{x^{j}}^{x^{l}}\partial_{x^{l}}J_{x^{i}}^{x^{k}}\frac{\partial\phi}{\partial x^{k}}+J_{x^{j}}^{x^{l}}J_{x^{i}}^{x^{k}}\frac{\partial^{2}\phi}{\partial x^{l}\partial x^{k}}\\
 & +J_{x^{j}}^{x^{l}}\partial_{x^{l}}J_{x^{i}}^{y^{k}}\frac{\partial\phi}{\partial y^{k}}+J_{x^{j}}^{x^{l}}J_{x^{i}}^{y^{k}}\frac{\partial^{2}\phi}{\partial x^{l}\partial y^{k}}+J_{x^{j}}^{x^{l}}\partial_{x^{l}}J_{x^{i}}^{t}\frac{\partial\phi}{\partial t}+J_{x^{j}}^{x^{l}}J_{x^{i}}^{t}\frac{\partial^{2}\phi}{\partial x^{l}\partial t}\\
 & +iJ_{x^{j}}^{y^{l}}\frac{\partial^{2}\phi}{\partial y^{l}\partial x^{i}}+J_{x^{j}}^{y^{l}}\partial_{y^{l}}J_{x^{i}}^{x^{k}}\frac{\partial\phi}{\partial x^{k}}+J_{x^{j}}^{y^{l}}J_{x^{i}}^{x^{k}}\frac{\partial^{2}\phi}{\partial y^{l}\partial x^{k}}+J_{x^{j}}^{y^{l}}\partial_{y^{l}}J_{x^{i}}^{y^{k}}\frac{\partial\phi}{\partial y^{k}}\\
 & +J_{x^{j}}^{y^{l}}J_{x^{i}}^{y^{k}}\frac{\partial^{2}\phi}{\partial y^{l}\partial y^{k}}+J_{x^{j}}^{y^{l}}\partial_{y^{l}}J_{x^{i}}^{t}\frac{\partial\phi}{\partial t}+J_{x^{j}}^{y^{l}}J_{x^{i}}^{t}\frac{\partial^{2}\phi}{\partial y^{l}\partial t}\\
 & +iJ_{x^{j}}^{t}\frac{\partial^{2}\phi}{\partial t\partial x^{i}}+J_{x^{j}}^{t}\partial_{t}J_{x^{i}}^{x^{k}}\frac{\partial\phi}{\partial x^{k}}+J_{x^{j}}^{t}J_{x^{i}}^{x^{k}}\frac{\partial^{2}\phi}{\partial t\partial x^{k}}\\
 & +J_{x^{j}}^{t}\partial_{t}J_{x^{i}}^{y^{k}}\frac{\partial\phi}{\partial y^{k}}+J_{x^{j}}^{t}J_{x^{i}}^{y^{k}}\frac{\partial^{2}\phi}{\partial t\partial y^{k}}+J_{x^{j}}^{t}\partial_{t}J_{x^{i}}^{t}\frac{\partial\phi}{\partial t}.
\end{align*}
In local coordinates we have $\Omega=hdz^{1}\wedge...\wedge dz^{n}$
for a holomorphic function $h$ and Equation (\ref{eq:AP1}) is equivalent
to 
\[
\det\left(\frac{\partial^{2}\phi}{\partial\overline{z}^{j}\partial z^{i}}\right)=|h|^{2}.
\]
Consider now the function 
\begin{align*}
G: & =\det\left(\frac{\partial^{2}\phi}{\partial\overline{z}^{j}\partial z^{i}}\right)-|h|^{2}\\
 & =\det\left(\frac{1}{4}\left(F_{i,j}+J_{t}^{x^{i}}J_{t}^{x^{j}}\frac{\partial^{2}\phi}{\partial t^{2}}\right)\right)-|h|^{2}.
\end{align*}
Obviously, $G$ is a polynomial of order $n$ in $\partial^{2}\phi/\partial t^{2}$
with functions defined on $W\times\mathbb{R}$ as coefficients. In
general, $G$ can be regarded as a function 
\[
G=G(p,t,\phi_{x},\phi_{y},\phi_{t},\phi_{xx},\phi_{yy},\phi_{xy},\phi_{tx},\phi_{ty},\phi_{tt}),
\]
where $\left(p,t\right)\in W\times\mathbb{R}$. By (\ref{eq:complex-structure-in-p})
and (\ref{eq:complex-structure-in-p-1}), there holds 
\begin{equation}
J|_{T_{p_{0}}W\times\mathbb{R}}=J_{\textnormal{st}}=\left(\begin{array}{cc}
0 & -\mathbb{I}\\
\mathbb{I} & 0
\end{array}\right),\label{eq:AP2}
\end{equation}
and we deduce that $J_{t}^{x^{i}}=0$ for $i=1,...,n-1$, and $J_{t}^{x^{n}}=1$.
Therefore, the equation $G=0$, restricted on $L$, can be solved
for $\partial^{2}\phi/\partial t^{2}$ because of 
\[
\left(D_{\phi_{tt}}G\right)|_{L}=\det\left(\begin{array}{ccc}
g_{11} & ... & g_{1(n-1)}\\
\vdots & \ddots & \vdots\\
g_{(n-1)1} & ... & g_{(n-1)(n-1)}
\end{array}\right)\not=0,
\]
where $g_{ij}$ are the metric coefficients of the metric $g$ on
$L$. Here we have used the initial condition $\phi|_{L}=\rho|_{L}$
and the fact that the Kähler metric induced by $(i/2)\partial\overline{\partial}\rho$
restricts to $g$ on $L$. Thus, by the real-analytic implicit function
theorem \cite{key-9}, the equation $G=0$ can be solved locally for
$\partial^{2}\phi/\partial t^{2}$, i.e. there exists an analytic
function $H=H(p,t,\phi_{x},\phi_{y},\phi_{xx},\phi_{yy},\phi_{xy},\phi_{t},\phi_{xt},\phi_{yt})$,
such that $\partial^{2}\phi/\partial t^{2}=H$ and $G(p,t,\phi_{x},\phi_{y},\phi_{xx},\phi_{yy},\phi_{xy},\phi_{t},\phi_{xt},\phi_{yt},H)=0$
in a neighbourhood of $L$ in $W\times\mathbb{R}$. As a result, the
initial-value problem
\begin{equation}
\begin{cases}
\det\left(\frac{\partial^{2}\phi}{\partial\overline{z}^{j}\partial z^{i}}\right) & =|h|^{2}\\
\phi(p,0) & =\rho(p,0)\\
\frac{\partial\phi}{\partial t}(p,0) & =\frac{\partial\rho}{\partial t}(p,0)
\end{cases}\label{eq:pde1}
\end{equation}
 is locally equivalent to
\begin{equation}
\begin{cases}
\frac{\partial^{2}\phi}{\partial t^{2}} & =H(p,t,\phi_{x},\phi_{y},\phi_{xx},\phi_{yy},\phi_{xy},\phi_{t},\phi_{xt},\phi_{yt})\\
\phi(p,0) & =\rho(p,0)\\
\frac{\partial\phi}{\partial t}(p,0) & =\frac{\partial\rho}{\partial t}(p,0).
\end{cases}\label{eq:pde2}
\end{equation}
Since all coefficient functions are real-analytic, by the Cauchy-Kovalewsky
theorem \cite{key-12} has a unique solution in a neighbourhood in
$W\times\mathbb{R}$ of each point $x\in L$. As $\rho$ is defined
globally and the local solutions of the Cauchy-Kovalewsky type are
unique, we obtain a solution $\phi$ of (\ref{eq:first-pde}) in a
neighbourhood of $L$ in $W\times\mathbb{R}$. Now, we check that
$\phi$ is strictly plurisubharmonic. Locally, we expand $\phi$ in
a power series in $t$, i.e. 
\[
\phi(x,y,t)=\phi_{0}(x,y)+\phi_{1}(x,y)t+\phi_{2}(x,y)t^{2}+....
\]
From $\phi(x,y,0)=\rho(x,y,0)$ and $\frac{\partial\phi}{\partial t}(x,y,0)=\frac{\partial\rho}{\partial t}(x,y,0)$,
we find that 
\[
\phi(x,y,t)=\rho(x,y,0)+\frac{\partial\rho}{\partial t}(x,y,0)t+t^{2}g(x,y,t),
\]
where $g(x,y,t)=\sum_{i=0}^{\infty}\phi_{i+2}(x,y)t^{i}$. Consequently,
we obtain 
\begin{align*}
\frac{\partial^{2}\phi}{\partial x^{j}\partial x^{k}}(x,0,0) & =\frac{\partial^{2}\rho}{\partial x^{j}\partial x^{k}}(x,0,0)\\
\frac{\partial^{2}\phi}{\partial x^{j}\partial y^{k}}(x,0,0) & =\frac{\partial^{2}\rho}{\partial x^{j}\partial y^{k}}(x,0,0)\\
\frac{\partial^{2}\phi}{\partial x^{j}\partial t}(x,0,0) & =\frac{\partial^{2}\rho}{\partial x^{j}\partial t}(x,0,0)\\
\frac{\partial^{2}\phi}{\partial y^{j}\partial t}(x,0,0) & =\frac{\partial^{2}\rho}{\partial y^{j}\partial t}(x,0,0)\\
\frac{\partial\phi}{\partial x^{j}}(x,0,0) & =\frac{\partial\rho}{\partial x^{j}}(x,0,0)\\
\frac{\partial\phi}{\partial y^{j}}(x,0,0) & =\frac{\partial\rho}{\partial y^{j}}(x,0,0)\\
\frac{\partial\phi}{\partial t}(x,0,0) & =\frac{\partial\rho}{\partial t}(x,0,0)
\end{align*}
and 
\[
\frac{\partial^{2}\phi}{\partial t^{2}}(x,0,0)=2g(x,0,0).
\]
Accounting of (\ref{eq:AP2}), we end up with 
\[
\left(\frac{\partial^{2}\phi}{\partial\overline{z}^{j}\partial z^{i}}\right)|_{L}=\left(\begin{array}{cccc}
\frac{\partial^{2}\rho}{\partial z^{1}\partial\overline{z}^{1}} & ... & \frac{\partial^{2}\rho}{\partial z^{1}\partial\overline{z}^{n-1}} & \frac{\partial^{2}\rho}{\partial z^{1}\partial\overline{z}^{n}}\\
\vdots & \ddots & \vdots & \vdots\\
\frac{\partial^{2}\rho}{\partial z^{n-1}\partial\overline{z}^{1}} & ... & \frac{\partial^{2}\rho}{\partial z^{n-1}\partial\overline{z}^{n-1}} & \frac{\partial^{2}\rho}{\partial z^{n-1}\partial\overline{z}^{n}}\\
\frac{\partial^{2}\rho}{\partial z^{n}\partial\overline{z}^{1}} & ... & \frac{\partial^{2}\rho}{\partial z^{n}\partial\overline{z}^{n-1}} & 2g(x,0,0)
\end{array}\right)|_{L}.
\]
The function $\phi$ is strictly plurisubharmonic if the matrix on
the left-hand side of the above equation is positive definite. In
fact, it is enough to show that this matrix is positive definite on
$L$, since in this case it follows that it is positive definite in
a neighbourhood of $L$. As $\rho$ is strictly plurisubharmonic,
the first $n-1$ principal minors of the matrix on the right-hand
side of the above equation are positive and it remains to show that
\[
\det\left(\frac{\partial^{2}\phi}{\partial\overline{z}^{j}\partial z^{i}}\right)|_{L}=\det\left(\begin{array}{cccc}
\frac{\partial^{2}\rho}{\partial z^{1}\partial\overline{z}^{1}} & ... & \frac{\partial^{2}\rho}{\partial z^{1}\partial\overline{z}^{n-1}} & \frac{\partial^{2}\rho}{\partial z^{1}\partial\overline{z}^{n}}\\
\vdots & \ddots & \vdots & \vdots\\
\frac{\partial^{2}\rho}{\partial z^{n-1}\partial\overline{z}^{1}} & ... & \frac{\partial^{2}\rho}{\partial z^{n-1}\partial\overline{z}^{n-1}} & \frac{\partial^{2}\rho}{\partial z^{n-1}\partial\overline{z}^{n}}\\
\frac{\partial^{2}\rho}{\partial z^{n}\partial\overline{z}^{1}} & ... & \frac{\partial^{2}\rho}{\partial z^{n}\partial\overline{z}^{n-1}} & 2g(x,0,0)
\end{array}\right)|_{L}>0.
\]
However, as $\phi$ solves (\ref{eq:pde1}) and $|h|^{2}$ is positive
on $L$, the conclusion readily follows, so $\phi$ is strictly plurisubharmonic.
Finally, we prove that the metric obtained from the Kähler form $\left(i/2\right)\partial\overline{\partial}\phi$
restricts on $L$ to the Riemannian metric on $L$. Because of 
\[
\det\left(\frac{\partial^{2}\phi}{\partial\overline{z}^{j}\partial z^{i}}\right)|_{L}=\det\left(\begin{array}{cccc}
g_{11} & ... & g_{1(n-1)} & g_{nn}\\
\vdots & \ddots & \vdots & \vdots\\
g_{(n-1)1} & ... & g_{(n-1)(n-1)} & g_{(n-1)n}\\
g_{n1} & ... & g_{n(n-1)} & 2g(x,0,0)
\end{array}\right)|_{L}
\]
and 
\[
\det\left(\frac{\partial^{2}\phi}{\partial\overline{z}^{j}\partial z^{i}}\right)|_{L}=|h|^{2}|_{L}=\det(g_{ij})
\]
we obtain $2g(x,0,0)=g_{nn}(x)$.

\end{document}